\documentclass[graybox]{svmult}

\usepackage{type1cm}        
\usepackage{makeidx}         
\usepackage{graphicx}        
\usepackage{multicol}        
\usepackage[bottom]{footmisc}

\usepackage{newtxtext}       %
\usepackage[varvw]{newtxmath}       

\usepackage{natbib}

\makeindex             

\let\varepsilon\varepsilon
\let\lif\to

\newcounter{matchct}
\newlength{\matchheight}
\newlength{\matchinc}
\newcommand{\matchstick}[2][1]{%
\setcounter{matchct}{0}
\setlength{\matchheight}{#1 mm}
\raisebox{1 mm}{%
\whiledo {\value{matchct} < #2}%
{%
  \,\rule[-.5\matchheight]{.5pt}{\matchheight}%
  \stepcounter {matchct}%
  \addtolength{\matchheight}{-\matchinc}
}}}

\def\eprint#1{}

\begin{document}

\title*{Hilbert's Program and Infinity}
\author{Richard Zach\orcidID{0000-0003-1633-8324}}
\institute{Richard Zach \at University of Calgary, Department of Philosophy, 2500 University Drive NW, Calgary, AB~T2N~1N4, Canada, \email{rzach@ucalgary.ca}}
%
%
\maketitle

\abstract{The primary aim of Hilbert's proof theory was to establish
the consistency of classical mathematics using finitary means only.
Hilbert's strategy for doing this was to eliminate the infinite (in
the form of unbounded quantifiers) from formalized proofs using the
so-called epsilon substitution method. The result is a formal proof
which does not mention or appeal to infinite objects or
``concept-formations.'' However, as later developments showed, the
consistency proof itself lets the infinite back into proof theory,
through a back door, so to speak. The paper outlines the epsilon
substitution method as an example of how proof-theoretic constructions
``eliminate the infinite'' from formal proofs, and how they aim to
establish conservativity and consistency. The proof also requires an
argument that this proof theoretic construction always works. This
second argument, however, requires possibly infinitary reasoning at
the meta-level, using induction on ordinal notations.}

\section{Introduction}
\label{sec:intro}

One of the aims of Hilbert's program was the elimination of the
infinite from the foundations of mathematics. Importantly, Hilbert did
not aim to eliminate the infinite from mathematics itself, i.e., he
did not aim to eliminate the use of the infinite in mathematical
practice. In fact, the whole point of the program was the
justification of the use of the infinite in mathematical practice by
showing that for foundational purposes, it could be altogether
avoided. The strategy was twofold: First, mathematics would be
\emph{formalized}, i.e., mathematical language was regimented in an
artificial language like that of predicate logic or type theory. The
work of Frege, Whitehead and Russell, Weyl, and others had shown how
this could be accomplished at least in principle for large parts of
the mathematics in use at Hilbert's time. This formalization results
in a formal image of mathematical practice: mathematical statements
become formulas in the regimented language, and mathematical proofs
become derivations from a (small) set of axioms also formulated in the
language. The second step was to establish, using means that avoid the
infinite as well, that these formal derivations never lead to a
contradiction, e.g., to a theorem of the form $\varphi \land \lnot
\varphi$ or a theorem the negation of which is also derivable from the
axioms such as~$\lnot 0=0$.\footnote{For general surveys of Hilbert's
proof-theoretic program, see
\cite{Kreisel1958a,Detlefsen1986,Zach2007a,Zach2023a}.}

The details of the strategy required the development of a formal
system in which to formalize mathematics and the development of
acceptable methods to establish its consistency. Proof theorists have
developed many approaches to both since Hilbert's day, but
it is instructive to look back at Hilbert's own. It will allow us to
highlight how it was supposed to work, and also where it stumbles. At
a very high level, the approach can be described as follows. The
formal system contains vocabulary the direct interpretation of which
seems to require the infinite: e.g., it might contain symbols for
(infinite) sets or functions and allow unrestricted quantification
over infinite domains. These are the ``ideal'' elements of the theory,
which a consistency proof aims to remove. One can identify a subsystem
which does not contain such ``ideal'' elements: a ``finitary'' kernel
of vocabulary and axioms which can be interpreted directly without
requiring the infinite. One such finitary subsystem is basic
quantifier free arithmetic. It contains symbols for (finite) natural
numbers only and no unrestricted quantification. The consistency of
this subsystem can be directly established on finitary grounds. The
consistency proof of the full system consists in a method for
transforming any derivation of a formula in the language of the safe
subsystem into a derivation that lies entirely within the subsystem,
together with a proof that the transformation always terminates and
yields a correct derivation.

The end-formula of derivations transformed by the method can be actual
theorems, such as $(7+5)=12$, but its scope extends to any prima facie
possible end formula such as the contradictory~$\lnot\,0=0$. The
consistency of the full system is then established. For suppose it
could derive a contradiction such as~$\lnot\,0=0$. The derivation in
the full system, which might rely on infinitary vocabulary or
infinitary axioms, would be transformed by the method contained in the
consistency proof into a derivation in the finitary subsystem. (The
consistency proof itself would have guaranteed, by finitary means,
that this method always works.) But the finitary subsystem itself can
directly be seen to be consistent, so it can't possibly prove
something false such as~$\lnot\,0=0$.

In what follows we will consider Hilbert's own preferred approach for
how to carry out the details. This approach
used the $\varepsilon$-calculus as the underlying system and the
$\varepsilon$-substitution method as the method used to transform
derivations in the full system to derivations in a finitary subsystem.
This will shed light on how Hilbert himself envisaged his consistency
program should be carried out. In Sect.~\ref{sec:ea} we will describe
the systems Hilbert used: first-order arithmetic formulated in the
$\varepsilon$-calculus as the full infinitary system and the
$\varepsilon$-free fragment as the finitary safe subsystem. In
Sect.~\ref{sec:es} we discuss the $\varepsilon$-substitution method
used to transform derivations in the full system to derivations in the
finitary subsystem. In Sec.~\ref{sec:ind} we will look at various
methods of induction, which are needed in order to establish that the
method works.

\section{Epsilon terms as ideal elements}
\label{sec:ea}

The vocabulary of first-order arithmetic contains terms for natural
numbers and the purely logical vocabulary. For definiteness, let us
consider a language with a symbol for $0$, for the successor
function~$'$, and for addition $+$ and multiplication~$\times$. Terms
are inductively defined as usual: $0$ is a term; if $t$ is a term, so
is $t'$; and if $t_1$ and $t_2$ are terms, so are $(t_1 + t_2)$ and
$(t_1 \times t_2)$. Thus, the following are terms:
\[0 \quad 0''' \quad (0''+0''') \quad (0'' \times (0''+0''')) \quad
\dots\] The terms that do not contain $+$ and $\times$ are special:
they are of the form $0^{\prime\cdots\prime}$, i.e., a $0$ followed by
some number of successor function symbols. We abbreviate a $0$
followed by $n$ such $'$~symbols as $\overline n$ and call it the
\emph{standard numeral} for~$n$. (If $n = 0$, then $\overline n$,
i.e., $\overline 0$, is just the symbol~$0$ by itself.) The standard
numerals are the simplest way to denote all the natural numbers in our
system.

The identity symbol~$=$ and the less than sign~$<$ allow us to form
the simplest possible formulas (we call them \emph{atomic}). If $t_1$
and $t_2$ are terms, then $t_1 = t_2$ and $t_1 < t_2$ are atomic
formulas. If $t_1$ and $t_2$ are numerals, then we can immediately
determine if such atomic formulas are intuitively true:
$0^{\prime\cdots\prime} = 0^{\prime\cdots\prime}$ is true if and only
if the number of primes on the left and right are the same, and
$0^{\prime\cdots\prime} < 0^{\prime\cdots\prime}$ if and only if there
are fewer primes on the left than on the right.

With the standard logical operators $\lnot$ (not), $\lor$ (or),
$\land$ (and), and $\to$ (if-then) we can form more complicated
formulas: if $A$ is a formula, so is $\lnot A$; if $A_1$ and $A_2$ are
formulas so are $(A_1 \lor A_2)$, $(A_1 \land A_2)$, and $(A_1 \to
A_2)$. If we only start from atomic formulas containing no $+$ or
$\times$, then using the usual truth tables, we can also directly
ascertain if such formulas are true. E.g., $(0' = 0' \to \lnot(0''' =
0'' \lor 0' < 0))$ is true since $0''' = 0''$ and $0' < 0$ are false,
which makes $(0''' = 0'' \lor 0' < 0)$ false, which makes $\lnot(0'''
= 0'' \lor 0' < 0)$ true, and a conditional with a true consequent is
true.

The addition and multiplication functions add a slight level of
complexity. But we know how to add and multiply (something we learned
in grade school) so every formula like $(\overline{5} + \overline{7}) =
(\overline{2} \times \overline{6})$ can be evaluated as true or false.
The formulas built from atomic formulas involving only $0$, $'$, $+$,
$\times$ and using only $=$, $<$, $\lnot$, $\lor$, $\land$, $\to$ can
all directly be seen to be true or false, using basic arithmetical
operations and the truth tables for the connectives. They are
finitary, ``real'' statements.

The ideal elements of first-order arithmetic are quantifiers. Hilbert
proposed a different system which replaces quantifiers by special
terms. In addition to the constant symbol~$0$ we allow terms to
contain variables ($x$, $y$, \dots). If $A(x)$ is a formula containing
the free variable~$x$, then $\varepsilon_x\,A(x)$ is an \emph{epsilon
term}. Intuitively, $\varepsilon_x\,A(x)$ stands for an arbitrary object
which satisfies~$A(x)$, if such an object exists; if no such object
exists then $\varepsilon_x\,A(x)$ may denote anything at all. E.g.,
$\varepsilon_x\,x<\overline{3}$ can be $0$, $1$, or $2$;
$\varepsilon_x\,x<x$ can be any number (since no number is less than
itself).

How can epsilon terms replace quantifiers? By stipulation,
$\varepsilon_x\,A(x)$ is a \emph{witness} for $A(x)$ (if one exists).
So clearly, $\exists x\,A(x)$ is true if and only if
$A(\varepsilon_x\, A(x))$ is. Also, $\varepsilon_x\,\lnot A(x)$ is a
\emph{counterexample} to~$A(x)$, if one exists. In that case, i.e., if
$\forall x\,A(x)$ is false, also $A(\varepsilon_x\,\lnot A(x))$ is
false. If no counterexample exists, then $\varepsilon_x\,\lnot A(x)$
must be a witness for $A(x)$, since $A(x)$ is true for all~$x$. So,
$\forall x\,A(x)$ is true if and only if $A(\varepsilon_x\,\lnot
A(x))$ is true.

Axiom systems for first-order logic with quantifiers require axioms
and/or special inference rules for $\forall$ and $\exists$. In the
\emph{epsilon calculus}, only one axiom schema is necessary to deal with
$\varepsilon$ terms:  \[A(t) \to A(\varepsilon_x\,A(x))\] Axioms of this
form are called \emph{critical formulas}. The only other axioms are
the usual ones for the propositional connectives, schemas for~$=$, and
the inference rule modus ponens: $A, A \to B \vdash B$.\footnote{See
\cite{AvigadZach2002} for a survey of the $\varepsilon$-calculus.}

To provide an axiom system for \emph{arithmetic}, we just need a few
axioms that say things about numbers, such as the following (we ignore
$<$ for simplicity):
\begin{align}
 \lnot\, 0 & = t'\\
 t_1' & = t_2' \to t_1 = t_2\\
 (t + 0) & = t\\
 (t_1 + t_2') & = (t_1 + t_2)'\\
 (t \times 0) & = 0\\
 (t_1 \times t_2') & = (t_1 + (t_1 \times t_2))
\end{align}
These schemas can all be written without any quantifiers or epsilon
terms, although $t$, $t_1$, and $t_2$ may contain epsilon terms.
However, if they don't, then we can see that every instance must be
true. E.g., whatever epsilon-free term $t$ may be, its value is some
natural number~$n$, so the value of $t'$ must be $n+1$, i.e., it must
be greater than~$0$. Thus, $0 = t'$ must be false, and hence $\lnot 0
= t'$ true.

The only axiom schema of usual first-order arithmetic that contains
problematic quantifiers is \emph{induction}:
\[
{[}A(0) \land \forall x(A(x) \lif A(x')){]} \lif \forall x\,A(x).
\]
Induction is equivalent to the \emph{least number principle}:
\[\exists x\,A(x) \lif \exists x(A(x) \land \forall y < x\,\lnot
A(y)).\] In arithmetic formulated in the epsilon calculus, we can add
an axiom to require that $\varepsilon_x\,A(x)$ is the \emph{least} witness
(if one exists):
\[A(t) \lif \varepsilon_x\,A(x) < t'.\] Adding all formulas of this form
completes the infinitary system of Hilbert's arithmetic. He called
them ``critical formulas of the second kind.'' The system has modus
ponens as the only inference rule. Only the critical formulas of first
and second kind mention epsilon terms, i.e., they are the only axioms
that mention ``ideal elements.'' All others, at least when they do not
contain epsilon terms, can intuitively be seen to be true. The
finitary, ``real'' subsystem thus excludes the critical formulas and
only allows epsilon-free instances of the remaining axiom schemas.

\section{Consistency proofs using the epsilon substitution method}
\label{sec:es}

The consistency of this real subsystem can be established easily. The
proof requires only induction on the complexity of terms and on the
length of derivations. One defines the \emph{value} of a term~$t$ not
containing variables or epsilons simply as the number it represents.
This is done inductively. The term $0$ has value zero. The value of
$t'$ is the result of adding~$1$ to the value of~$t$. The value of
$(t_1 + t_2)$ is the sum of the values of $t_1$ and $t_2$, and the
value of $(t_1 \times t_2)$ the product of their values. Next we can
define a truth value for atomic formulas containing only such terms.
A sentence of the form $t_1 = t_2$ is true if the values of $t_1$ and
$t_2$ are the same, and false otherwise. A sentence of the form $t_1 <
t_2$ is true if the value of $t_1$ is less than~$t_2$. Complex
formulas containing $\lnot$, $\lor$, $\land$, and $\to$ get truth
values assigned inductively as well, using the usual truth tables.
Since real formulas involve no quantifiers, variables, or epsilons,
this gives us a finitary notion of truth for all real formulas. To
determine if a formula is true or false we just need to be able to
count, add and subtract natural numbers, and apply the finite truth
tables of $\lnot$, $\lor$, $\land$, and~$\to$. Truth and falsity of
real formulas is mechanically decidable in finitely many steps.

To establish consistency we show that every formula derivable in the
real subsystem is true under this interpretation. We use proof by
induction on the length of derivations. Every derivation in the real
subsystem is a sequence of formulas each one of which is either an
axiom or follows from previous ones in the sequence by the one rule of
inference, modus ponens. To show consistency, we establish that every
end-formula of such a derivation is true. We do this by showing that
all axioms are true and that modus ponens leads from true formulas to
true formulas. Since we're assuming the derivation does not involve
critical formulas, we just have to observe that every propositional
axiom is a tautology (and hence true regardless of the truth values of
formulas it is made up of). For instance, say we have an axiom of the form $(A
\land B) \to A$. Either $A$ is true or false. If it is true, the axiom
is a conditional with a true consequent; if it is false, $(A \land
B)$ is false and so the axiom is a conditional with a false
antecedent. In either case, the axiom is true. If the axiom is an
instance of one of the arithmetical axioms we can likewise prove by
finitary reasoning about the values of the terms $t_1$ and $t_2$ that
the values of left and right sides of each equation must always be the
same.

This approach does not directly extend to terms, formulas, and
derivations involving epsilon terms, simply because epsilon terms
cannot be assigned values in a finitary way. Hilbert proposed to deal
with this problem is as follows. We can restrict the question of
showing that the full, ideal system is consistent to the question of
showing that, say, $\lnot 0=0$ cannot have a derivation. The real system
already proves $0=0$ (it is an instance of an identity axiom), so if
it also proved $\lnot 0=0$ it would be inconsistent. And any
inconsistent system proves any formula whatsoever, including $\lnot
0=0$. Hilbert's idea is to show that if the ideal system proves $\lnot
0=0$, then so does the real system. But $\lnot 0=0$ is a false real
formula, so can't have a purely real proof by the consistency proof of
the real subsystem.

Hilbert's approach is easy to grasp in the case where our purported
ideal derivation of $\lnot 0=0$ contains a single epsilon term, say,
$\varepsilon_x A(x)$.\footnote{Hilbert first articulated this approach
in \cite{Hilbert1923}. It was then mainly carried out by his student
\cite{Ackermann1924,Ackermann1940}. A textbook presentation of
Ackermann's proof is given in \cite{HilbertBernays1970}. A modern
presentation and analysis of Ackermann's proof is given by
\cite{Moser2006}. \cite{Tait1965} gives an alternative
$\varepsilon$-substitution method. The history of consistency proofs
using the $\varepsilon$-calculus in the 1920s is discussed in
\citet{Zach2003,Zach2004a}.} Then all the axioms that are not also
axioms of the real subsystem are just critical formulas involving
$\varepsilon_x A(x)$, i.e., of the form of either 
\[A(t_1) \to A(\varepsilon_x A(x)) \qquad\text{or}\qquad A(s_1) \lif
\varepsilon_x\,A(x) < s_1'.\] Take the purported derivation of $\lnot
0=0$ and replace every occurrence of $\varepsilon_x A(x)$ in it
by~$0$. This uniform replacement turns every propositional, identity,
and arithmetical axiom into another instance of the same axiom, and so
these become true formulas under the finitary interpretation of terms
and formulas not involving epsilons. The critical formulas turn into
formulas of the form of either
\[A(t_2) \to A(0) \qquad\text{or}\qquad A(s_2) \lif 0 < s_2'.\] (Note
that the terms $t_1$ and~$s_1$ may turn into different terms~$t_2$, $s_2$
after the replacement because $\varepsilon_x A(x)$ may occur in~$t_1$
or $s_1$
and be replaced by~$0$.)

The formulas on the right are clearly all true, since the value
of~$s_2'$ is always greater than~$0$. To determine if the formulas on
the left are true, we apply our finitary interpretation of real
formulas. If all of them are true, we are in the same situation as in
the consistency proof for real derivations. We have a derivation of
$\lnot 0=0$, a false formula, from the true axioms of the real
subsystem plus other formulas (the results of the replacement of
$\varepsilon_x A(x)$ in all critical formulas) which are also all true.
That can't be, so at least one of the formulas of the form $A(t_2) \to
A(0)$ must be false. Since a conditional is only false if the
antecedent is true and the consequent is false, we can conclude that
$A(t_2)$ is true (and $A(0)$ false). Let $k$ be the value of the
term~$t_2$. By evaluating $A(\overline 1)$, $A(\overline 2)$, \dots,
$A(\overline k)$, find the least~$n$ such that $A(\overline n)$ is
true. Take the original ideal derivation containing $\varepsilon_x A(x)$,
and now replace $\varepsilon_x A(x)$ everywhere by~$\overline n$. This
again results in a derivation from the (true) axioms of the real
subsystem plus formulas of the form
\[A(t_3) \to A(\overline n) \qquad\text{or}\qquad A(s_3) \lif \overline n < s_3'.\]
By the way we have found~$n$, we know that $A(\overline n)$ and hence
any of the formulas of the first kind are true.
Since $n$ is least among the numbers for which~$A(x)$ is true,
whatever the value~$\ell$ of~$s_3$ is, if $A(\overline \ell)$ is true,
$n < \ell+1$, and so the formulas of the second kind are also true.

Before we go on to discuss the general case, a few observations.
First, both this proof and the preceding proof of the consistency of
the real subsystem are constructive. This is obscured by the fact that
they begin from the assumption that we have a derivation of a
contradiction such as $\lnot 0=0$, and so they have a whiff of
``indirectness'' about them. But the statement we are proving is a
negative statement: there can be no derivation of~$\lnot 0=0$.
Constructive proofs of negative statements simply have this form:
assume there were such a derivation, and reduce this assumption to a
contradiction. Both claims and proofs can easily be reformulated as
positive claims which have consistency as a consequence. In the
case of the real subsystem, the positive claim is: every derivation in
the real subsystem contains only true formulas (``true'' understood in
the finitary meaning we have provided). It entails consistency
because it shows that $\lnot 0=0$ cannot be the end-formula of a real
derivation. In the case of ideal derivations containing only a single
$\varepsilon$-term, the positive claim is: for every derivation of a real
formula using only critical formulas belonging to a single
$\varepsilon$-term, there is a real derivation of the same formula. The
same argument works to establish this result. Start with a derivation
of a real formula~$F$ containing only critical formulas belonging to
$\varepsilon_x A(x)$. Replace $\varepsilon_x A(x)$ everywhere by~$0$ and
verify if all critical formulas have been replaced by true formulas.
If so, we have a derivation of~$F$ from true real formulas. If not,
$A(0)$ is false but at least one $A(t_2)$ is true: find the least~$n$
such that $A(\overline n)$ is true and replace $\varepsilon_x A(x)$
instead by~$n$. The result is also guaranteed to be a derivation
of~$F$ from true real formulas.

The only step missing now is to show that the true formulas $A(t_2) \to
A(0)$ (or $A(t_3) \to A(\overline n)$) have real derivations themselves.
This is not hard to do. The key step is to prove, by induction on the
complexity of terms~$t$, that the arithmetic axioms of the real
subsystems provide derivations of $t = \overline k$, where $k$ is the
value of~$t$. Together we have shown something stronger than
consistency: we have shown conservativity of the ideal system over its
real subsystem, i.e., that if $F$ is a real formula and has a
derivation in the ideal system (only involving a single
$\varepsilon$-term), then $F$ has a derivation already in the real
subsystem. Consistency follows, since the real subsystem only derives
true formulas and $\lnot 0=0$ is not true.

So far we have only sketched how to establish conservativity of the
ideal system over the real for derivations involving a single epsilon
term. Of course, ideal derivations can contain any number of epsilon
terms. Importantly, these epsilon terms can themselves contain other
epsilon terms, and the way in which they can be nested inside one
another can be arbitrarily complicated. To illustrate, let's consider
the translation of a simple formula containing two existential
quantifiers into the $\varepsilon$-calculus. Recall that the translation
of $\exists x\,A(x)$ is $A(\varepsilon_x\,A(x))$. We apply this
translation successively to obtain a translation of $\exists
x\exists y\,A(x,y)$. First we focus on the $\exists x$ quantifier:
\begin{align*}
\fbox{$\exists x$}\,\exists y\, & A(\fbox{$x$},y)\\
\intertext{and replace it with an
$\varepsilon$-term. We obtain}
\exists y\,& A(\fbox{$\varepsilon_x\, \exists z\,A(x,z)$}, y).\\
\intertext{Note that we've replaced the $\exists y$ inside the $\varepsilon$-term by $\exists z$ to avoid confusion with the $\exists y$ remaining outside. In the second step, we replace the new quantifier $\exists z$:}
\exists y\,& A(\varepsilon_x\, \fbox{$\exists z$}\,A(x,\fbox{$z$}), y)\\
\intertext{thus turns into}
\exists y\, &A (\varepsilon_x\, A(x, \fbox{$\varepsilon_z\,A(x,z)$}), y).\\
\intertext{Finally, we replace $\exists y$ in}
\fbox{$\exists y$}\, & A(\varepsilon_x\, A(x, \varepsilon_z\,A(x,z)), \fbox{$y$})\\
\intertext{by the corresponding $\varepsilon$-term to obtain the result:}
& A(\varepsilon_x\, A(x,\varepsilon_z\,A(x,z)), \fbox{$\varepsilon_y\,A(\varepsilon_x\, A(x,\varepsilon_z\,A(x,z)), y)$}).
\end{align*}
As you can imagine, derivations in the ideal system containing the equivalent
in the language of the $\varepsilon$-calculus may contain very complex
$\varepsilon$-terms. E.g., the equivalent of $\exists x\exists z\, x<y$
would read:
\[
\varepsilon_x\, \left[x < \varepsilon_z \left(x < z\right)\right] < \varepsilon_y\, \left[\varepsilon_x\,
\left(x < \varepsilon_z\left(x <z\right)\right) < y\right].
\]
More importantly for our purposes, Hilbert's approach for proving
conservativity becomes vastly more complicated. This is not just
because the number of $\varepsilon$-terms in critical formulas can be
arbitrary in the general case, but these $\varepsilon$-terms can also be
contained in each other in complex ways, and these containment
relationships affect how we must go about finding replacements of $\varepsilon$-terms
by numerals so that all critical formulas become true real formulas.
Hilbert's approach proceeds by starting with a guess for a numerical
value for each $\varepsilon$-term occurring in the derivation, and then
changing that guess if the guess doesn't pan out, i.e., turn every
critical formula into a true formula (which then is also derivable in
the real subsystem). If only one $\varepsilon$-term occurs, the initial
guess is~$0$, which may get changed to the least~$n$ such that
$A(\overline n)$ is true if $A(0)$ is false.

Here's a first complication if more than one $\varepsilon$-term occurs.
Suppose there are two $\varepsilon$-terms, $\varepsilon_x\, A(x)$ and
$\varepsilon_y B(y)$, with just one critical formula each:
\[
A(t_1) \to A(\varepsilon_x\, A(x)) \qquad\text{and}\qquad B(s_1) \to B(\varepsilon_y B(y)).
\]
Say we replace both by~$0$ as before to get
\[
A(t_2) \to A(0) \qquad\text{and}\qquad B(s_2) \to B(0).
\]
Suppose also that $A(t_2)$ is true and $A(0)$, $B(s_2)$,
and~$B(0)$ are false, making the first critical formula false and the
second true. Following our previous strategy, we'd ``correct'' our
guess for $\varepsilon_x\, A(x)$ to the least~$n$ such that $A(\overline
n)$ is true and redo the replacement. This would yield
\[
A(t_3) \to A(\overline n) \qquad\text{and}\qquad B(s_3) \to B(0)
\]
(we still replace $\varepsilon_y\,B(y)$ by~$0$ since the second critical
formula was true already). The first critical formula is now true.
Note that $t_1$ and $s_1$ may contain both $\varepsilon_x\,A(x)$ and
$\varepsilon_y\,B(y)$. So the results of replacing those two
$\varepsilon$-terms in~$s_1$, resulting in terms $s_2$ and~$s_3$, may be
different. In one case we replace $\varepsilon_x\,A(x)$ by~$0$ and
in the other by~$\overline n$. Consequently, even if $B(s_2)$ is false,
$B(s_3)$~may not be. So our improved guess made the first critical
formula true, but perhaps the second critical formula false. In this
simple case, we could now just find the least~$m$ such
that~$B(\overline m)$ is true, and get two true formulas
\[
A(t_4) \to A(\overline n) \qquad B(s_4) \to B(\overline m).
\]
Even this simple case shows that there may be many steps involved
until we arrive at a ``solving substitution'' of numerals for
$\varepsilon$-terms.

But this is not the only, or even the most important, difficulty.
Recall that the system includes not just axiom schemas for the propositional
connectives, but also for identity. One such axiom schema is $t_1=t_2
\to s(t_1) = s(t_2)$, i.e., terms of equal value may be substituted
for one another in a term~$s$ while preserving its value. Now recall
that $\varepsilon$-terms may occur as subterms in each other. E.g., in
the example above (the translation of $\exists
x\exists y\,A(x,y)$), we ended up with a formula containing the two $\varepsilon$-terms
\begin{align*}
& \varepsilon_x\, A(x,\varepsilon_z\,A(x,z)) \text{ and}\\
& \varepsilon_y\, A(\fbox{$\varepsilon_x\, A(x,\varepsilon_z\,A(x,z))$}, y),
\end{align*}
where the first is contained in the second. It quickly becomes clear
that taking our guesses to apply to full $\varepsilon$-terms without
taking into account which $\varepsilon$-terms are contained in which
others might wreak havoc, and we have to be more careful. E.g., let's
say for simplicity we have three $\varepsilon$-terms
\[
\varepsilon_x\, A(x), \qquad \varepsilon_y\,B(y, \overline n), \qquad\text{and}\qquad\varepsilon_y\,B(y, \varepsilon_x\, A(x)),
\]
and we simply guessed $0$ as the value for all three.
Then suppose our derivation contains an instance of an equality axiom, say,
\[\overline n=\varepsilon_x\,A(x) \to \varepsilon_y\,B(y, \overline n) = \varepsilon_y\,B(y,
\varepsilon_x\, A(x)).\]
After replacing $\varepsilon$-terms with our guesses, this turns into
$\overline n = 0 \to 0=0$. It has a false antecedent, so we don't have
to worry about it. But let's say our guess of~$0$ for $\varepsilon_x\, A(x)$ is
incorrect, and we determine that $n$ happens to be least so that
$A(\overline n)$ holds. Our new replacement turns the identity axiom
into $\overline n = \overline n \to 0 = 0$. This is still true. Let's
say in the next step we determine that the replacement for
$\varepsilon_y\,B(y, \overline n)$ can't be $0$ and must be $m > 0$. Now we'd
get the false $\overline n = \overline n \to \overline m = \overline
0$. This shows that we can't simply treat the guesses for
$\varepsilon_x\, A(x)$, $\varepsilon_y\,B(y, \overline n)$ and $\varepsilon_y\,B(y,
\varepsilon_x\, A(x))$ as independent. Instead, we must always replace
$\varepsilon_y\,B(y, \overline n)$ and $\varepsilon_y\,B(y,
\varepsilon_x\, A(x))$ by the same numeral if we replace $\varepsilon_x\,
A(x)$ by~$\overline n$. The obvious way to do this is to not guess
numerals for whole $\varepsilon$-terms occurring in the proof, but to
guess numerals for any potential $\varepsilon$-term resulting from a
replacement of a nested term by a numeral. In our case, we should not
guess a value for $\varepsilon_y\,B(y,
\varepsilon_x\, A(x))$ as a whole, but for $\varepsilon_y\,B(y,
0)$ if our guess for $\varepsilon_x\, A(x)$ is~$0$, and for $\varepsilon_y\,B(y,
\overline n)$ if that guess is~$n$. Of course, it's not always just a
single value for which we need a guess. In addition to $\varepsilon_y\,B(y,
\varepsilon_x\, A(x))$, our derivation might also contain $\varepsilon_y\,B(y,
\varepsilon_z\, C(z))$, and if our guess for $\varepsilon_z\, C(z)$ is $m$
we also need a guess for $\varepsilon_y\,B(y,
\overline m)$.

It would take us too far to describe how Ackermann's epsilon
substitution method deals with all of these issues technically.
Suffice it to say that it is a complicated procedure, which requires
keeping track of substitutions for not just the $\varepsilon$-terms actually
occurring in the derivation, but also for new ones that result from
those by substituting nested $\varepsilon$-terms by numerals. In very
broad strokes, the procedure actually deals with what are called
$\varepsilon$-\emph{types}. An $\varepsilon$-type is a
$\varepsilon$-term which contains no nested subterms except free
variables, and in which no free variable occurs more than once. An
$\varepsilon$-substitution assigns a numerical function to each
$\varepsilon$-type which has the property that its value is non-zero
only for finitely many (tuples of) arguments. Interpreting an
$\varepsilon$-type as such a function on the one hand allows us to
determine a numerical value for every $\varepsilon$-term in the
derivation, and on the other hand solves the difficulties noted above.
The procedure begins by starting with an $\varepsilon$-substitution
that assigns functions that are $0$ for very possible argument to all
applicable $\varepsilon$-types, and improves this in each step by
adding new non-zero values to some of them, but also by sometimes
starting over with the default $0$-valued function for some
$\varepsilon$-types. Note that crucially an $\varepsilon$-substitution
so defined is finite: it consists of finitely many
$\varepsilon$-types, with finitely many tuples of arguments and their
non-zero values. For all other (infinitely many) arguments, the
function has value zero; these do not have to be kept track of.

The number of steps the procedure has to run through before a solving
substitution is found far exceeds the number of $\varepsilon$-terms
actually occurring in the derivation, since each step may generate a
new tuple of arguments for which a non-zero value is determined. Since
some such new guesses require previous guesses to be abandoned, the
procedure backtracks a lot, and care must be taken that it doesn't go
in circles. The important point is that this can be done. There is a
deterministic procedure which always produces a solving substitution
after finitely many steps, and which transforms a derivation of any real
formula involving $\varepsilon$-terms into one without, i.e., a real
derivation. The procedure moreover manipulates only finitary objects,
i.e., finite sequences of symbols (e.g., derivations and mappings from
a finite number of $\varepsilon$-types to finitely many
argument--value pairs).

\section{Consistency proofs, induction, and infinity}
\label{sec:ind}

We now have a relatively clear idea of how finitary consistency proofs
work. Typically, they consist of a procedure to transform a derivation
in the system to be proved consistent (the ideal system, e.g., formal
arithmetic with $\varepsilon$-terms and critical formulas) into a
derivation of the same end-formula in a ``safe'' system (e.g., the
real subsystem of arithmetic without $\varepsilon$-terms or quantifiers).
Of course, the end-formula must be a formula allowed in the safe
subsystem, e.g., a real formula involving just arithmetical
operations, equality, and truth-functional operators. Such a formula
has a finitary meaning and can be finitarily evaluated as true or
false, and we can prove, finitarily, that the safe subsystem only
produces true formulas of this kind. In particular, it cannot prove
$\lnot 0=0$, and so is consistent. By removing $\varepsilon$-terms from
derivations, such a procedure eliminates those expressions of the
language that only have an infinitary meaning, i.e., involve infinite
search or quantification. In that sense, a consistency proof is a
procedure for removing the infinite from derivations in the ideal
system.

There is, however, a place where the infinite creeps back in. That is
in the verification that the procedure always works. We don't just
want a procedure that happens to work in every case (and which is
itself finitary, operates on finite expressions, etc.) but we also
want to \emph{verify} that it works. Such a verification will have to
use meta-theoretical inference principles, such as induction. Suppose
$P$ is some property~$P$ of natural numbers. Regular induction on
natural numbers can be formulated as an inference of the following
form. If, for all natural numbers~$n$, $P(n)$ holds if $P(k)$ holds
for all $k < n$, then $P(n)$ holds for all~$n$. To prove that $P(n)$
holds for all~$n$, we show that, for arbitrary~$n$, $P(n)$ holds \emph{on
the assumption that it holds for all~$k<n$}. Since consistency proofs
aren't directly about numbers but about other things such as terms or
derivations, to use this principle we have to associate numbers as
``measures'' or ``weights'' with the expressions about which we are
proving things.

Let's consider a simple example. Part of the verification that the
real subsystem is consistent required that this system proves every
true real formula. For such a proof we'd have to verify that a number
of things are true, such as: If $n$ is the value of~$t$, then the
formula $t = \overline n$ has a derivation, or: if $n \neq m$, then
the formula $\lnot\, \overline n = \overline m$ has a derivation. Let's
see how an inductive proof of the latter claim would go. To make things
a bit simpler, assume $n<m$. Then we can write $m = n + \ell + 1$ where
$\ell\ge 0$. Let's prove that, for all $n$, $\lnot\, \overline n =
\overline{n + \ell + 1}$ has a derivation. The ``measure'' we assign here
to the claim ``$\lnot\, \overline n = \overline{n + \ell + 1}$ has a
derivation'' is just~$n$ (and so is independent of~$\ell$). To prove it
by induction we assume that the claim holds for all $k<n$, and then
use this assumption (the ``inductive hypothesis'') in the proof that
it holds for~$n$ as well. 

So let $n$ be an arbitrary natural number. There are two cases. The
number $n$ might be $0$, or it might be greater than~$0$. In the first
case where $n=0$, the term $\overline{n+\ell+1}$ is
just~$\overline{\ell+1}$. By definition, $\overline{\ell+1}$ identical to the
term~$\overline{\ell}'$. Thus, the formula $\lnot\, \overline n =
\overline{n + \ell + 1}$ is identical to the formula $\lnot\, 0 =
\overline{\ell}'$, which is one of the axiom schemas, so it has a
derivation.

Now consider the second case: $n>0$. Then $n=k+1$ where $k = n-1$.
Clearly, $k<n$. The inductive hypothesis says that $\lnot\,
\overline{k} = \overline{k+\ell+1}$ has a derivation. The following
are axioms:
\begin{align*}
\overline{k}' = \overline{k+\ell+1}' & \to \overline{k} = \overline{k+\ell+1}\\
(\overline{k}' = \overline{k+\ell+1}' & \to \overline{k} = \overline{k+\ell+1}) \to (\lnot\,\overline{k} = \overline{k+\ell+1} \to \lnot\,\overline{k}' = \overline{k+\ell+1}')
\end{align*}
(The second is an instance of contraposition, i.e., $(A \to B) \to
(\lnot B \to \lnot A)$, a tautology.) By applying modus ponens twice,
we obtain a derivation of $\lnot\,\overline{k}' = \overline{k+\ell+1}'$,
which is identical to $\lnot\,\overline n = \overline{n + \ell + 1}$.

Why does induction work? It is supposed to establish $P(n)$ for an
arbitrary~$n$. It does so by proving that $P(n)$ holds if $P(n_1)$
holds for some $n_1 < n$, which in turn follows from the hypothesis
that $P(n_2)$ holds for some $n_2 < n_1$, etc. (In the example, $n_1 =
n-1$, $n_2 = n-2$, etc.) Eventually, this sequence of claims for ever
decreasing measures must end, since a decreasing sequence of natural
numbers cannot decrease forever. (In our case, it ends when we reach
$P(0)$ which is proved without the use of the induction hypothesis but
rather follows from the fact that $\lnot\, 0 = t'$ is an axiom.) 

The natural numbers are, as we say, \emph{well-ordered:} there are no
indefinitely decreasing sequences of natural numbers. Induction on
natural numbers works because the natural numbers (in their natural
ordering~$<$) are well-ordered. We can tell that the natural numbers
are well-ordered because we know that there are only $n$ natural
numbers less than~$n$, so a decreasing sequence of natural numbers
starting with $n$ can contain at most $n+1$ elements before it must
reach~$0$ (namely when the decrease in each step is minimal: $n$,
$n-1$, $n-2$, $n-3$, \dots, $0$). Induction on the natural numbers is
clearly finitarily justified, because the well-ordering of the natural
numbers is immediately evident, given that the natural numbers just
are those that can be reached from $0$ by a finite number of
applications of the ``$+1$'' successor operation.

Before we go on to consider more complex and
perhaps no longer finitarily justifiable induction principles and
well-orders, let's pause to consider how our explicitly inductive
proof that $\lnot\, \overline n = \overline m$ has a real derivation
compares to the transformation procedure in the
$\varepsilon$-substitution method. After all, the latter was not
formulated as an inductive proof, but rather sketched as a procedure
that successively transforms derivations into others. Even supposing
it works as advertised, what does it have to do with an induction?

Let's consider our simple example again, but reformulate it as a
transformation procedure. The transformation procedure works on
partial derivations of $\lnot\, \overline{n} = \overline{n+\ell+1}$.
They are partial in the sense that not every formula in them has to be
an axiom: One single formula in them can be of this form without
itself being justified by modus ponens from preceding formulas. In
fact, the formula $\lnot\, \overline{n} = \overline{n+\ell+1}$ by
itself is such a partial derivation, and would be the partial
derivation we would start from. The transformation procedure is the
following. Suppose you have a partial derivation of the described
form, and let $\lnot\, \overline{n} = \overline{n+\ell+1}$ be the
single formula in it that is not justified by modus ponens. If this
formula is such that $n=0$, then we already have a complete derivation
of its end-formula, since then the formula in question is really
$\lnot\, 0 = \overline{\ell}'$, which is an axiom. Otherwise, $n =
k+1$. Write down a new partial derivation by taking the old one and
add the following four formulas at the beginning:
\begin{align*}
  \lnot\, &
\overline{k} = \overline{k+\ell+1}  \tag{1}\\
  & \overline{k}' = \overline{k+\ell+1}' \to \overline{k} = \overline{k+\ell+1} \tag{2}\\
   (&\overline{k}' = \overline{k+\ell+1}'  \to \overline{k} = \overline{k+\ell+1}) \to (\lnot\,\overline{k} = \overline{k+\ell+1} \to \lnot\,\overline{k}' = \overline{k+\ell+1}') \tag{3}\\
  \lnot\, & \overline{k} = \overline{k+\ell+1} \to \lnot\,\overline{k}' = \overline{k+\ell+1}' \tag{4}
\end{align*}
where $k=n-1$ as before. Formulas 2 and 3 are axioms. Formula 4 is
justified by modus ponens from 2 and 3, and the previously unjustified
$\lnot\, \overline{n} = \overline{n+\ell+1}$, i.e.,
$\lnot\,\overline{k}' = \overline{k+\ell+1}'$ is justified by modus
ponens from 1 and 4. The new partial derivation has a single
unjustified formula, namely formula~1.

We have described a procedure which transforms a partial derivation
into a new partial derivation with a new, single unjustified formula.
The only case where it does not result in a new partial derivation is
the case where the unjustified formula is of the form $\lnot 0 =
\overline{\ell}'$, which is an axiom: so when the procedure halts we
have obtained a complete derivation where every formula is an axiom or
justified by modus ponens. \emph{If we knew} that this procedure always
terminates, we would know that it would produce, starting with the
target formula $\lnot\, \overline{n} = \overline{n+\ell+1}$ on its
own, a complete correct derivation thereof. How do we know that it
always terminates? By stipulation and inspection of the procedure, all
partial derivations involved contain only a single line of the form
$\lnot\, \overline{n} = \overline{n+\ell+1}$ that is not justified.
Let's say that this $n$ is the ``weight'' of such a partial
derivation. The procedure replaces a partial derivation of weight~$n$
by one with weight~$k = n-1$, i.e., a smaller weight. The sequence of
partial derivations produced by the procedure corresponds to a
decreasing sequence of their weights. Since these weights are natural
numbers and those are well-ordered, there can be no infinite
decreasing sequence of weights, and consequently no infinite sequence
of partial derivations they correspond to. Hence, the procedure must
eventually stop and arrive at a complete derivation.

The lesson for the $\varepsilon$-substitution method is that the
specification of the procedure by itself does not yet constitute a
consistency (or more precisely, conservativity) proof. It also
requires a proof that it always terminates. Such a proof can be given
if we can assign weights to the derivations produced such that each
step produces a new derivation with smaller weight. The procedure only
halts when all $\varepsilon$-terms and critical formulas have been
removed from the proof and replaced by numerals, i.e., when we have
reached a real derivation of the same end-formula. Unfortunately, the
``weights'' required are much more complex than single numbers. For
instance, we cannot simply assign as weights the number of
$\varepsilon$-terms or the number of critical formulas since these numbers may well
\emph{increase} in the course of finding the solving $\varepsilon$-substitution.

What kind of measure could we even consider as appropriate weights? One
still relatively simple measure would be pairs $\langle n_1,
  n_2\rangle$ of natural numbers with the following ordering:
  $\langle n_1, n_2 \rangle < \langle k_1, k_2\rangle$ if and only if
  either \begin{enumerate}
  \item $n_1 = k_1$ but $n_2 < k_2$, or 
  \item $n_1 < k_1$ 
  \end{enumerate}
For instance, in this ordering $\langle 1, 2\rangle < \langle 1,
6\rangle$ by~(1), and $\langle 1, 6 \rangle < \langle 2, 1\rangle$
by~(2). It is not hard to see that pairs of natural numbers are
well-ordered by this ordering. For suppose a decreasing sequence
starts with $\langle n_1, n_2\rangle$. By using clause~(1), we can
pass to a lesser pair by decreasing the second component and leaving
the first the same, i.e., by passing to $\langle n_1, k_2\rangle$ with
$k_2 < n_2$. But there are only finitely many numbers $< n_2$.
Eventually we must use clause~(2) and decrease the first component,
i.e., pass to $\langle k_1, k_2\rangle$ with $k_1<n_1$. Importantly,
when this happens there is no restriction on~$k_2$; it may be any
natural number whatsoever and need not be $< n_2$. Whatever $k_2$ is,
there are only finitely many numbers less than it, so eventually we
must use clause~(2) again and decrease the first component. But $n_1$
is also a natural number, so we can decrease the first component also
only finitely many times.

This ordering of pairs would be used, for instance, if we wanted to
prove that $t = \overline{\ell}$ has a real derivation if $t$ is a term
containing $+$ (but no~$\times$) and $\ell$ is the value of~$t$. Our
proof would proceed as before. We start with a partial derivation,
e.g., $t = \overline \ell$ itself. Step by step, add formulas to the
beginning of this partial derivation, until we obtain a complete
derivation in which every formula is either an axiom or justified by
modus ponens. Each partial derivation contains a single unjustified
line of the form $t = \overline{\ell}$ (the term~$t$ changes, but the
value~$\ell$ stays the same). Each such derivation can be
assigned a weight of the form $\langle n_1, n_2\rangle$. Without going
into the details, we might pick as $n_1$ the number of occurrences of
$+$ in~$t$, and as $n_2$ the sum of values of all
numerals~$\overline{k}$ occurring as the right argument to~$+$. For
instance, the weight of $(\overline{3} + (\overline{10} +
\overline{1})) = \overline{14}$ would be $\langle 2, 1\rangle$. Each
step would produce a new derivation in which either the first
component stays the same, but the second component is reduced (e.g.,
the unjustified formula might now be $(\overline{3} + (\overline{10} +
0)') = \overline{14}$ with weight $\langle 2, 0\rangle$) or the first
component is reduced (in the following step we would have $(\overline{3} +
\overline{11}) = \overline{14}$ with weight $\langle 1, 11\rangle$).

The set of pairs of numbers with their associated ordering as weights
is more complex than the natural numbers with their natural ordering.
The objects themselves are still finite (just pairs of numbers). The
ordering itself is finitary (i.e., it is decidable whether $\langle
n_1, n_2 \rangle < \langle k_1, k_2\rangle$ without any unbounded
search or unrestricted quantification). It can be seen that the
ordering is a well-order (the first component can only be decreased
finitely many times, and each time it happens, the second component
can only be decreased finitely many times). All this together
justifies thinking of induction on pairs of numbers as finitary, and
the procedure described that turns $t = \overline{\ell}$ into a real
derivation thereof as a finitary procedure which moreover can be seen,
finitarily, to always terminate. But induction on pairs of numbers is
not complicated enough to show that the $\varepsilon$-substitution
procedure always terminates.

To get a sense of how complex the well-ordering required to show
termination of $\varepsilon$-substituion is, it helps to consider a
standard way of measuring the complexity of well-orders. Such a way is
provided by the set-theoretic idea of ordinals. Ordinals are one of
two ways in which infinite numbers appear in set theory, the other
being cardinals. While cardinals measure the size of sets (including
infinite ones), ordinals measure the length of well-orders.
Set-theorists define arithmetical operations on ordinals. Among
others, we have a successor function (``$+1$''), addition,
multiplication, and exponentiation. Any well-order is measured by an
ordinal. The simplest well-order is a set containing a single element.
Any such well-order has ordinal~$1$. If any well-order has ordinal
$\alpha$, then the well-order with a single element added ``at the
end'' has the ordinal $\alpha+1$. For instance, well-orders with 1, 2,
3, \dots{} elements, e.g.,
\[a_1 \qquad a_1 < a_2 \qquad a_1 < a_2 < a_3 \qquad \dots\]
have ordinals 1, 2, 3, \dots{}

The smallest infinite well-order is that of the natural numbers, and
its ordinal is~$\omega$. Let's represent such a well-order graphically as
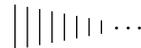
\[\matchstick[6]{8}\cdots\]
``Matchsticks'' further left represent elements of the well-order that
are smaller than those to the right, even though the length of the
matchsticks decrease as the elements themselves increase in the order.
Addition of ordinals corresponds to stringing well-orders together,
e.g., $\omega+\omega$ might be represented as
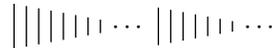
\[
  \matchstick[6]{8}\cdots \ \matchstick[5]{7}\cdots
\]
Multiplication of ordinals, e.g., $\alpha \cdot \beta$ corresponds to
replacing every element of a well-order measured by $\beta$ with an
entire well-order corresponding to $\alpha$. For instance, to get
$\omega\cdot 2$, take two elements, say, $0 < 1$, and replace each
element with a copy of the natural numbers. We'll get the same
well order as displayed above---unsurprisingly, as we'd expect $\omega
+ \omega$ and $\omega \cdot 2$ to be the same.\footnote{Note that
addition and multiplication of ordinals is not commutative. For
instance, $2 \cdot \omega$ corresponds to the well-order resulting
from replacing every natural number ($\omega$) with two copies where
the first is less than the second. But the well-order $0_1 < 0_2 < 1_1
< 1_2 < 2_1 < 2_2 < \dots$ is not relevantly different from the
natural numbers themselves, and also has the ordinal~$\omega$.} A
well-order with ordinal $\omega \cdot \omega$ (or $\omega^2$) would
then be represented as
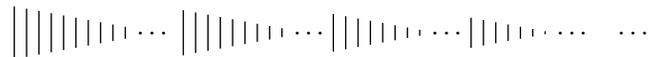
\[
\matchstick[7]{10}\cdots \ \matchstick[6]{9}\cdots \matchstick[5]{8}\cdots \matchstick[4]{7}\cdots \quad\cdots
\]

We've now come far enough to understand that the ordinal associated
with the well-order of pairs of numbers defined above is $\omega^2$.
Each pair of numbers $\langle n_1, n_2\rangle$ has a place in our
matchstick diagram: the first component $n_1$ tells you which of the
$\omega$-blocks the pair is in, the second component tells you how far
into that $\omega$-block it is located.

Set-theoretically, ordinals are certain kinds of infinite sets.
Assuming set-theory, we could similarly justify our procedure above by
assigning ordinals themselves as weights, namely the ordinal
$\omega\cdot n_1 + n_2$ to any partial derivation to which we've
previously assigned the pair $\langle n_1, n_2\rangle$. All of these
ordinals are less than $\omega^2$. Since $\omega^2$, like any ordinal,
is itself a well-ordered set (in fact, it just is the set of all
ordinals less than it, set-theoretically speaking), we could then also
conclude that our procedure always terminates. This would not be a
finitary justification, since it would make use of infinite objects
(set-theoretic ordinals) and a fact justified by infinitary set theory
(all ordinals are well-ordered). To see \emph{finitarily} that our
procedure terminates we have to give a finitary justification as we
did above. As a notational convenience, proof theorists often talk as
if they are doing this, i.e., use ordinals as weights and proceed by
induction ``up to $\alpha$.'' But what they really mean (or at least
should mean) is that they are talking about \emph{notations} for
ordinals, such as pairs of numbers like $\langle 3, 5\rangle$, or
strings of symbols like ``$\omega\cdot 3 + 5$'' (and finitarily
acceptable order relations on them like the one for pairs we described
above).

Our procedure above has left out the multiplication symbol~$\times$
from the term~$t$. Of course, it is also true that $t = \overline{\ell}$
has a real derivation whenever $\ell$ is the value of~$t$, even if $t$
contains both $+$ and~$\times$. The proof of this, i.e., the
corresponding procedure that eventually produces such a derivation,
would be more complicated. We would require a more complicated measure
to prove that it terminates. We could use triples $\langle n_1, n_2,
n_3\rangle$ of numbers ordered as follows:
$\langle n_1, n_2, n_3\rangle < \langle k_1, k_2, k_3\rangle$ if and
only if either $n_1 < k_1$ (in the ordinary sense of~$<$) or $n_1 =
k_1$ and $\langle n_2, n_3\rangle < \langle k_2, k_3\rangle$ (in the
sense we've given to $<$ for pairs). This would still be finitary. We
could convince ourselves of this as follows: If a decreasing sequence
of triples starts with $\langle n_1, n_2, n_3\rangle$, there are only
finitely many times the first component can decrease (since $n_1$ is
finite). Each time this happens, the second and third component can
become arbitrarily large. But as we know, pairs of numbers are
well-ordered, so for each of the finitely many times the first
component decreases, there are also only finitely many times the pairs
making up the second and third component can decrease.

The ordinal corresponding to triples of numbers is~$\omega^3$. We can
keep developing more complicated orderings in this way, e.g.,
quadruples of numbers correspond to $\omega^4$ and so on. As we did
for triples, we can ``bootstrap'' our way up this hierarchy of more
and more complex well-orders and produce finitary arguments that these
measures are well-ordered. But we can't bootstrap the same way to the
next level of complexity,~$\omega^\omega$. We can still produce finitary objects and a
finitary ordering on them which corresponds to the
ordinal~$\omega^\omega$. (The measures in this case could be sequences
of natural numbers of arbitrary but finite length.) 

Even this complicated measure is not enough to prove that the
$\varepsilon$-substitution method terminates. What would be required
is a set of measures, or ``ordinal notations'' as proof theorists call
them, that correspond to the ordinal~$\varepsilon_0$. Recall that the
$\varepsilon$-substitution method proceeds by successively improving
``guesses'' of certain functions for all $\varepsilon$-``types''
relevant to the derivation. Each such guess consists of a finite
sequence of $\varepsilon$-types and a finite sequence of
argument--value pairs for each of them. It is these collective guesses
that are each assigned an ordinal notation. The proof that the method
terminates consists in showing that each step in the procedure
produces a new $\varepsilon$-substitution with an ordinal notation
that is less than the preceding one. If we know that these ordinal
notations are well-ordered, we know that the procedure cannot run
forever.

The ordinal $\varepsilon_0$ is the least ordinal greater than all of
$\omega$, $\omega^\omega$, $\omega^{\omega^\omega}$, \dots{} Again, it
is not terribly hard to describe finite objects with a decidable order
on them that has this ordinal. But what is hard, and perhaps
impossible, is to give an argument that this ordering is a well-order
which withstands finitary scrutiny. Which ordinal notations can be
finitarily justified has been the topic of some debate in the
literature. For instance, \cite{Tait1981} has given an influential
argument that the limit of finitary reasoning is characterized
by~$\omega^\omega$. On the other hand, \citet[\S11]{Takeuti1975} has
attempted a finitary justification for induction up
to~$\varepsilon_0$.\footnote{See \cite{Zach2003} on early uses of
ordinal notations exceeding $\omega^\omega$ in \cite{Ackermann1924}.
See Ch.~8 of \cite{MancosuGalvanZach2021} for an in-depth discussion
of finitary well-orders and ordinal notations up to~$\varepsilon_0$.
\cite{DarnellThomas-Bolduc2022a} discuss Takeuti's
justification of~$\varepsilon_0$.}

\section{Conclusion}

Proof-theoretic results about the consistency of formal systems of
ideal, infinitary mathematics usually proceed by transforming
derivations in such systems into derivations in other formal systems
which are less problematic. ``Less problematic'' means that the target
systems can be seen independently to be consistent. Such
transformations typically establish more than mere consistency: they
establish conservativity for formulas in the real system. For
instance, they transform ideal derivations of real formulas into real
derivations of the same formulas. Moreover, it can be shown, by
finitary means, that the real system derives only true real formulas,
i.e., is \emph{sound}. Since contradictions are not true, the
soundness of the real system entails its consistency. If the
transformation can also be seen to always produce a real derivation,
i.e., always terminates, this then establishes the consistency of the
ideal system.

The transformation involved in this proof-theoretic argument removes
finitarily problematic elements from the proof. This is especially
clear in the case of the $\varepsilon$-substitution method, where the
only questionable elements are $\varepsilon$-terms. The substitution
method replaces these by numerals. For other systems and proof
transformations, the matter is less transparent. For instance, in the
case of standard Peano arithmetic~$\mathbf{PA}$, the usual proofs,
going back to Gentzen, remove applications of the induction rule as
well as quantifiers from derivations of a real formula.\footnote{A
detailed proof of Gentzen's approach can be found in
\cite{MancosuGalvanZach2021}.}

All such transformations are very complicated, and the proofs that
they terminate require ``weights,'' i.e., ordinal notations, that
correspond to quite large ordinals. In the case of systems of
first-order arithmetic, this ordinal is~$\varepsilon_0$. This is not a
mere accident or due to a lack of better proof ideas. \cite{Gentzen1943}
showed how ordinal notations can themselves be formalized in
arithmetic, i.e., there is a sentence~$I(\varepsilon_0)$ of
first-order arithmetic that ``says that'' ordinal notations up to
$\varepsilon_0$ are well-ordered.\footnote{The matter is of
formalizing ordinal notations is subtle; the results here mentioned do
depend on the definition of the notations in arithmetic and not just
on the ordinal itself as \cite{Kreisel1976} has highlighted. Textbook
presentations of the results mentioned here can be found, e.g., in
\cite{Takeuti1975}.} He also showed that Peano arithmetic derives the
claim $I(\varepsilon_0) \lif \mathrm{Con}(\mathbf{PA})$, i.e., that
``induction up to $\varepsilon_0$'' implies the consistency of
\textbf{PA}. Since $\mathbf{PA}$ does not prove its own consistency,
by G\"odel's second incompleteness theorem, it can't prove
$I(\varepsilon_0)$. At the same time, $\mathbf{PA}$ does prove
$I(\alpha)$ for $\alpha<\varepsilon_0$ (as long as the ordinal
notations for~$\alpha$ are represented in the same direct way). Again,
by the second incompleteness theorem, no smaller~$\alpha$ suffices to
prove the consistency of~$\mathbf{PA}$. If it did, such an argument
would be expected to also be formalizable in~$\mathbf{PA}$, and so
then $\mathbf{PA}$ would derive both $I(\alpha)$ and $I(\alpha) \lif
\mathrm{Con}(\mathbf{PA})$ contradicting the second incompleteness
theorem.

The proof-theoretic methods as well as the ordinal notation systems
used for them are true feats of mathematical ingenuity. Nevertheless, the
question of the extent to which they allow philosophical conclusions
about the security of ideal mathematics remains open to debate. We
have seen how they remove the ``ideal,'' i.e., (expressions for)
infinite objects from ideal derivations. However, they do not
eliminate all traces of infinity, or infinitary reasoning. The place
where infinitary reasoning definitely still plays a role is in the
justification of induction along the ordinal notations used, i.e., that
the weights assigned to derivations are well-ordered. The problem here
is not that these weights themselves or their ordering is finitarily
suspect; they are not. What is questionable is that finitary methods
alone suffice to establish that they are well-ordered.

\bigskip

\noindent\textbf{Acknowledgments.} Thanks to Paolo Mancosu and Koray
Akçagüner for comments on a draft of this paper. The material was
presented at the ``Infini mathématique'' seminar at the University of
Paris~I, at the Logic Café, University of Vienna, the \emph{Gödel and
Kant on Mathematics and Physics} workshop at the University of
Tübingen, and the \emph{Logic in Philosophy: Incompleteness
and Intuition} workshop at the University of Rijeka. I am grateful to
the organizers and audiences for their questions and comments.

\bibliographystyle{spbasic}
\bibliography{infinite}
\end{document}